\documentclass[times, authoryear, 12pt]{elsarticle}

\usepackage{geometry}
\usepackage[labelfont=bf]{caption}
\usepackage{amsmath,amsfonts,amssymb,amsthm,color,epsfig,graphicx,url}
\usepackage{appendix}
\usepackage[utf8]{inputenc}
\usepackage[english]{babel}

\urldef\myurl\url{www.ualberta.ca/~dwiens/}
\newtheorem{theorem}{Theorem}

\newtheorem{lemma}{Lemma}

\makeatletter
\def\ps@pprintTitle{  \let\@oddhead\@empty
  \let\@evenhead\@empty
  \def\@oddfoot{\reset@font\hfil\thepage\hfil}
  \let\@evenfoot\@oddfoot
}
\makeatother
\journal{U. Alberta preprint series}
\begin{document}

\date{\today}

\begin{frontmatter}
\title{Minimum Variance Designs \linebreak With Constrained Maximum Bias}

\author[A1]{Douglas P. Wiens\corref{mycorrespondingauthor}}
\address[A1]{Mathematical \& Statistical Sciences,
	University of Alberta,
	Edmonton, Canada,  T6G 2G1}

\cortext[mycorrespondingauthor]{E-mail: \url{doug.wiens@ualberta.ca}. Supplementary material is at {\myurl}.}

\begin{abstract}
Designs which are minimax in the presence of
model misspecifications have been constructed so as to minimize the maximum,
over classes of alternate response models, of the integrated mean squared
error of the predicted values. This mean squared error decomposes into a
term arising solely from variation, and a bias term arising from the model
errors. Here we consider the problem of designing so as to minimize the variance of the
predictors, subject to a bound on the maximum (over model misspecifications)
bias. We consider as well designing so as to minimize the maximum bias,
subject to a bound on the variance. We show that solutions to both problems
are given by the minimax designs, with appropriately chosen values of their
tuning constants. Conversely, any minimax design solves each problem for an appropriate choice of the bound on the maximum bias or on the variance. 
\end{abstract}

\end{frontmatter}

\section{Introduction and summary \label{section: intro}}

An experimental design to fit a particular model is\textit{\ robust} if its
performance is stable under perturbations of the true model. The
theory of robustness of design was largely initiated by \cite{bd59}, who
investigated the robustness of some classical experimental designs in the
presence of certain model inadequacies, e.g.\ designs optimal for a low
order polynomial response when the true response was a polynomial of higher
order. \cite{h75} derived designs for straight line regression, robust in
the presence of alternate response functions. \cite{w90a,w92} extended these
results to multiple regression responses and in a variety of other
directions -- see \cite{w15} for a summary of these and other approaches to
robustness of design.

Designs which are \textit{minimax} in the presence of model
misspecifications aim to minimize the maximum, over classes of alternate
response models, of the integrated mean squared error (\textsc{imse}) of the
predicted values. In this note we first briefly review the theory of such
designs, and discuss a decomposition of the \textsc{imse} into a convex
combination of two terms -- one arising solely from variation, and the other
arising from the bias due to the model errors. It is often the case that
designs which yield a small value of one of these terms do poorly with
respect to the other. For instance the designs optimal with respect to the
variance-based alphabetic optimality criteria concentrate their mass at a
minimal number of points, and thus fare poorly, with large biases, in the
presence of model misspecifications. On the other hand uniform-like designs
reduce the bias while increasing the variance. Thus in \S \ref{section:
bounded} we propose two associated problems: (i) design so as to minimize
the integrated variance of the predictors, subject to a bound on the maximum
(over model misspecifications) bias, and (ii) design so as to minimize the
maximum bias, subject to a bound on the variance. We show that solutions to
both problems are given by the minimax designs, for appropriately chosen
values of the mixing parameter. Conversely, any minimax design solves both
problems for appropriate choices of the bounds on the maximum bias or
variance.

Examples of these ideas as they apply to polynomial regression, and methods of implementation, are studied in \S \ref{section:
examples}. The \textsc{matlab} code used to prepare these examples is
available on the author's personal website.

\section{Minimax robustness of design \label{section: minimax}}

The general minimax design problem is phrased in terms of an approximate
regression response 
\begin{equation}
E\left[ Y\left( \boldsymbol{x}\right) \right] \approx \boldsymbol{f}^{\prime
}\left( \boldsymbol{x}\right) \boldsymbol{\theta },  \label{approx}
\end{equation}%
for $p$ regressors $\boldsymbol{f}$, each functions of $q$ independent
variables $\boldsymbol{x}$, and a parameter $\boldsymbol{\theta }$. Since (%
\ref{approx}) is an approximation the interpretation of $\boldsymbol{\theta }
$ is unclear; we \textit{define }this target parameter by 
\begin{equation}
\boldsymbol{\theta }=\arg \min_{\boldsymbol{\eta }}\int_{\mathcal{X}}\left( E%
\left[ Y\left( \boldsymbol{x}\right) \right] -\boldsymbol{f}^{\prime }\left( 
\boldsymbol{x}\right) \boldsymbol{\eta }\right) ^{2}\mu \left( d\boldsymbol{x%
}\right) ,  \label{theta def}
\end{equation}%
where $\mu \left( d\boldsymbol{x}\right) $ represents either Lebesgue
measure or counting measure, depending upon the nature of the \textit{design
space} $\mathcal{X}$ with $\int_{\mathcal{X}}\mu \left( d\boldsymbol{x}%
\right) <\infty $. We then define $\psi \left( \boldsymbol{x}\right) =E\left[
Y\left( \boldsymbol{x}\right) \right] -\boldsymbol{f}^{\prime }\left( 
\boldsymbol{x}\right) \boldsymbol{\theta }$. This results in the class of
responses $E\left[ Y\left( \boldsymbol{x}\right) \right] =\boldsymbol{f}%
^{\prime }\left( \boldsymbol{x}\right) \boldsymbol{\theta }+\psi \left( 
\boldsymbol{x}\right) $, with -- by virtue of (\ref{theta def}) -- $\psi $
satisfying the orthogonality requirement 
\begin{equation}
\int_{\mathcal{X}}\boldsymbol{f}\left( \boldsymbol{x}\right) \psi \left( 
\boldsymbol{x}\right) \mu \left( d\boldsymbol{x}\right) =\boldsymbol{0}.\ 
\label{orthogonality}
\end{equation}%
Assuming that $\mathcal{X}$ is rich enough that the matrix\ $\boldsymbol{A}%
=\int_{\mathcal{X}}\boldsymbol{f}\left( \boldsymbol{x}\right) \boldsymbol{f}%
^{\prime }\left( \boldsymbol{x}\right) \mu \left( d\boldsymbol{x}\right) $
is invertible, the parameter defined by (\ref{theta def}) and (\ref%
{orthogonality}) is unique.

We identify a design with its design measure -- a probability measure $\xi
\left( d\boldsymbol{x}\right) $ on $\mathcal{X}$. Define 
\begin{equation*}
\boldsymbol{M}_{\xi }=\int_{\mathcal{X}}\boldsymbol{f}\left( \boldsymbol{x}%
\right) \boldsymbol{f}^{\prime }\left( \boldsymbol{x}\right) \xi \left( d%
\boldsymbol{x}\right) ,\text{ \ }\boldsymbol{b}_{\psi ,\xi }=\int_{\mathcal{X%
}}\boldsymbol{f}\left( \boldsymbol{x}\right) \psi \left( \boldsymbol{x}%
\right) \xi \left( d\boldsymbol{x}\right) ,
\end{equation*}%
and assume $\xi $ is such that $\boldsymbol{M}_{\xi }$ is invertible. For an
anticipated design of $n$, not necessarily distinct, points the covariance
matrix of the least squares estimator $\boldsymbol{\hat{\theta}}$, assuming
homoscedastic errors with variance $\sigma _{\varepsilon }^{2}$, is $\left(
\sigma _{\varepsilon }^{2}/n\right) \boldsymbol{M}_{\xi }^{-1}$, and the
bias is $E[ \boldsymbol{\hat{\theta}-\theta }] =\boldsymbol{M}%
_{\xi }^{-1}\boldsymbol{b}_{\psi ,\xi }$; together these yield the mean
squared error (\textit{mse}) matrix 
\begin{equation*}
\text{\textsc{mse}}\left[ \boldsymbol{\hat{\theta}}\right] =\frac{\sigma
_{\varepsilon }^{2}}{n}\boldsymbol{M}_{\xi }^{-1}+\boldsymbol{M}_{\xi }^{-1}%
\boldsymbol{b}_{\psi ,\xi }\boldsymbol{b}_{\psi ,\xi }^{\prime }\boldsymbol{M%
}_{\xi }^{-1}
\end{equation*}%
of the parameter estimates, whence the \textit{mse} of the predicted values $%
\hat{Y}\left( \boldsymbol{x}\right) =\boldsymbol{f}^{\prime }\left( 
\boldsymbol{x}\right) \boldsymbol{\hat{\theta}}$ is%
\begin{equation*}
\text{\textsc{mse}}\left[ \hat{Y}\left( \boldsymbol{x}\right) \right] =\frac{%
\sigma _{\varepsilon }^{2}}{n}\boldsymbol{f}^{\prime }\left( \boldsymbol{x}%
\right) \boldsymbol{M}_{\xi }^{-1}\boldsymbol{f}\left( \boldsymbol{x}\right)
+\left( \boldsymbol{f}^{\prime }\left( \boldsymbol{x}\right) \boldsymbol{M}%
_{\xi }^{-1}\boldsymbol{b}_{\psi ,\xi }\right) ^{2}+\psi ^{2}\left( 
\boldsymbol{x}\right) .
\end{equation*}%
\ 

\noindent A loss function that is commonly employed is the \textit{%
integrated mse} of the \textit{predictions}: 
\begin{eqnarray}
\text{\textsc{imse}}\left( \xi |\psi \right) &=&\int_{\mathcal{X}}\text{%
\textsc{mse}}\left[ \hat{Y}\left( \boldsymbol{x}\right) \right] \mu \left( d%
\boldsymbol{x}\right)  \notag \\
&=&\frac{\sigma _{\varepsilon }^{2}}{n}tr\left( \boldsymbol{AM}_{\xi
}^{-1}\right) +\boldsymbol{b}_{\psi ,\xi }^{\boldsymbol{\prime }}\boldsymbol{%
M}_{\xi }^{-1}\boldsymbol{AM}_{\xi }^{-1}\boldsymbol{b}_{\psi ,\xi }+\int_{%
\mathcal{X}}\psi ^{2}\left( \boldsymbol{x}\right) \mu \left( d\boldsymbol{x}%
\right) .  \label{imse}
\end{eqnarray}%
The dependence on $\psi $ is eliminated by adopting a \textit{minimax}
approach, according to which one first maximizes (\ref{imse}) over a
neighbourhood of the assumed response, i.e. over $\psi$. This neighbourhood is constrained by (%
\ref{orthogonality}) and by a bound $\int_{\mathcal{X}}\psi ^{2}\left( 
\boldsymbol{x}\right) \mu \left( d\boldsymbol{x}\right) \leq \tau ^{2}/n$,
required so that errors due to bias and to variation remain of the same
order, asymptotically. Define $\psi _{0}\left( \boldsymbol{x}\right) =\sqrt{n%
}\psi \left( \boldsymbol{x}\right) /\tau $ and $\nu =\tau ^{2}/\left( \sigma
_{\varepsilon }^{2}+\tau ^{2}\right) $. Then $\max_{\psi }$\textsc{imse}$%
\left( \xi |\psi \right) $ is $\left( \sigma _{\varepsilon }^{2}+\tau
^{2}\right) /n$ times 
\begin{equation*}
I_{\nu }\left( \xi \right) =\left( 1-\nu \right) \text{\textsc{var}}\left(
\xi \right) +\nu \text{\textsc{maxbias}}\left( \xi \right) ,
\end{equation*}%
where \textsc{var}$\left( \xi \right) =tr\boldsymbol{AM}_{\xi }^{-1}$ is proportional to the
integrated variance of the predictors and \textsc{maxbias}$\left( \xi
\right) =\max_{\psi _{0}}\text{\textsc{bias}}\left( \xi |\psi _{0}\right) $,
where 
\begin{equation}
\text{\textsc{bias}}\left( \xi |\psi _{0}\right) =\boldsymbol{b}_{\psi
_{0},\xi }^{\prime }\boldsymbol{M}_{\xi }^{-1}\boldsymbol{AM}%
_{\xi }^{-1}\boldsymbol{b}_{\psi _{0},\xi }+1  \label{conditional bias}
\end{equation}%
is the integrated (squared) bias, with $\psi _{0}$ constrained by (\ref%
{orthogonality}) and by $\int_{\mathcal{X}}\psi _{0}^{2}\left( \boldsymbol{x}%
\right) \mu \left( d\boldsymbol{x}\right) =1$.

We note that our implicit assumption of i.i.d. errors is itself minimax -- 
\cite{w24a} has shown that, in large classes of alternate corrrelation structures constrained only by a bound on an induced matrix norm, such as the maximum eigenvalue, of the covariance matrix of the errors, the i.i.d. structure is in fact least favourable. Thus minimizing under this assumption is a minimax procedure.
 
\section{Robust Bounded Maximum Bias and Bounded Variance designs \label%
{section: bounded}}

Let $\Xi $ be a class of designs on $\chi $, for instance all probability
measures on $\left[ -1,1\right] $ -- requiring appropriate approximations to
make them implementable -- or all exact designs, i.e.\ with integer
allocations, on a finite design space. For%
\begin{equation}
\min_{\xi \in \Xi }\text{\textsc{maxbias}}\left( \xi \right) \leq b^{2}\leq
\max_{\xi \in \Xi }\text{\textsc{maxbias}}\left( \xi \right) ,
\label{b-range}
\end{equation}%
consider the problem

\begin{description}
\item[(B):] Minimize \textsc{var}$\left( \xi \right) $ in the class $\Xi
_{b}\subset $ $\Xi $ of designs for which \textsc{maxbias}$\left( \xi
\right) \leq b^{2}$.
\end{description}

\noindent We call a solution to (B) a \textit{Robust Bounded Bias} design
with bias bound $b^{2}$, denoted \textsc{rbb}($b^{2}$). For%
\begin{equation}
\min_{\xi \in \Xi }\text{\textsc{var}}\left( \xi \right) \leq s^{2}\leq
\max_{\xi \in \Xi }\text{\textsc{var}}\left( \xi \right) ,  \label{s-range}
\end{equation}%
consider the problem

\begin{description}
\item[(S):] Minimize \textsc{maxbias}$\left( \xi \right) $ in the class $\Xi
_{s}\subset $ $\Xi $ of designs for which \textsc{var}$\left( \xi \right)
\leq s^{2}$.
\end{description}

\noindent We call a solution to (S) a \textit{Robust Bounded Variance}
design with variance bound $s^{2}$, denoted \textsc{rbv}($s^{2}$).\medskip

For $\nu \in \left[ 0,1\right] $ define $\xi _{\nu }=\arg \min_{\xi }I_{\nu
}\left( \xi \right) $. Set $b^{2}\left( \nu \right) =$ \textsc{maxbias}$%
\left( \xi _{\nu }\right) $ and $s^{2}\left( \nu \right) =$ \textsc{var}$%
\left( \xi _{\nu }\right) $. Note that then%
\begin{equation*}
\xi _{\nu }\in \Xi _{b\left( \nu \right) }\cap \Xi _{s\left( \nu \right) }.
\end{equation*}%
Below we prove and discuss Theorem \ref{thm: solution}%
. This asserts that solutions to (B) and (S) are given by $\xi _{\nu }$, for
some $\nu \in \left[ 0,1\right] $, for any $b^{2}$ and any $s^{2}$.
Conversely, any $\xi _{\nu }$ is a solution to (B) and (S) for appropriate
values of $b^{2}$ and $s^{2}$.

\begin{theorem}
\label{thm: solution} (a) Robust Bounded Bias designs with bias bound $b^{2}$
satisfying (\ref{b-range}) are given by 
\begin{equation*}
\xi =\left\{ 
\begin{array}{cc}
\xi _{\nu }, & \text{if }b^{2}=b^{2}\left( \nu \right) \text{, for }0\leq
\nu \leq 1\text{;} \\ 
\xi _{0}, & \text{if }b^{2}\geq b^{2}\left( 0\right) \text{.}%
\end{array}%
\right.
\end{equation*}%
$\newline
$(b) Robust Bounded Variance designs with variance bound $s^{2}$ satisfying (%
\ref{s-range}) are given by 
\begin{equation*}
\xi =\left\{ 
\begin{array}{cc}
\xi _{\nu }, & \text{if }s^{2}=s^{2}\left( \nu \right) \text{, for }0\leq
\nu \leq 1\text{;} \\ 
\xi _{1}, & \text{if }s^{2}\geq s^{2}\left( 1\right) \text{.}%
\end{array}%
\right.
\end{equation*}
\end{theorem}

\noindent  \textbf{Proof of Theorem \protect\ref{thm: solution} } \medskip

By the \textit{I-optimal} design (\cite{s77}) we mean the minimizer of $%
I_{0}\left( \xi \right) $, i.e.\ of the Integrated Variance of the Predicted
Values. By the \textit{uniform} design we mean the design $\xi \left( d%
\boldsymbol{x}\right) \propto \mu \left( d\boldsymbol{x}\right) $. These
designs play special roles -- they turn out to be $\xi _{0}$ and $\xi _{1}$,
respectively. 

Theorem \protect\ref{thm: solution} is immediate from the following three lemmas.

\begin{lemma}
\label{lemma: special}The design $\xi _{0}$, minimizing $I_{0}\left( \xi
\right) =$ \textsc{var}$\left( \xi \right) $ in $\Xi $, is I-optimal and the
design $\xi _{1}$, minimizing $I_{1}\left( \xi \right) =$ \textsc{maxbias}$%
\left( \xi \right) $ in $\Xi $ is uniform.
\end{lemma}

\noindent \textbf{Proof}: That $\xi _{0}$ is the I-optimal design follows
from the definition: $I_{0}\left( \xi \right) =$ \textsc{var}$\left( \xi
\right) $. By (\ref{conditional bias}), \textsc{maxbias}$\left( \xi \right)
\geq 1$. This lower bound is attained by the uniform design, since then, by (\ref{orthogonality}),
\begin{equation*}
\boldsymbol{b}_{\psi ,\xi }=\int_{\mathcal{X}}\boldsymbol{f}\left( 
\boldsymbol{x}\right) \psi \left( \boldsymbol{x}\right) \xi \left( d%
\boldsymbol{x}\right) \propto \int_{\mathcal{X}}\boldsymbol{f}\left( 
\boldsymbol{x}\right) \psi \left( \boldsymbol{x}\right) \mu \left( d%
\boldsymbol{x}\right) =\boldsymbol{0}.
\end{equation*} 
 \hfill $\square $

By Lemma \ref{lemma: special}, the lower bounds of the ranges (\ref{b-range}%
) and (\ref{s-range}) are attained by $\xi _{1}$ and $\xi _{0}$ respectively.

\begin{lemma}
\label{lemma: BB}(a) $\xi _{\nu }$ is a solution to (B) for $b^{2}\left( \nu
\right) $. (b) If $b^{2}\geq b^{2}\left( 0\right) =$ \textsc{maxbias}$\left(
\xi _{0}\right) $ then $\xi _{0}$ is a solution to (B) for $b^{2}$.
\end{lemma}

\noindent \textbf{Proof}: (a) For any $\xi \in \Xi _{b\left( \nu \right) }$
we must have that \textsc{var}$\left( \xi \right) \geq $ \textsc{var}$\left(
\xi _{\nu }\right) $, since otherwise%
\begin{equation*}
I_{\nu }\left( \xi \right) =\left( 1-\nu \right) \text{\textsc{var}}\left(
\xi \right) +\nu \text{\textsc{maxbias}}\left( \xi \right) <\left( 1-\nu
\right) \text{\textsc{var}}\left( \xi _{\nu }\right) +\nu b^{2}\left( \nu
\right) =I_{\nu }\left( \xi _{\nu }\right) ,
\end{equation*}%
a contradiction. Thus $\xi _{\nu }\in \Xi _{b\left( \nu \right) }$ minimizes 
\textsc{var}$\left( \xi \right) $ in $\Xi _{b\left( \nu \right) }$, i.e.\ is
\textsc{rbb}($b^{2}\left( \nu \right) $). \newline
(b) For such $b^{2}$ we have that \textsc{maxbias}$\left( \xi _{0}\right)
\leq b^{2}$ so that $\xi _{0}\in \Xi _{b}$. As well, 
\begin{equation}
\text{\textsc{var}}\left( \xi _{0}\right) =\min_{\xi \in \Xi }\text{\textsc{%
var}}\left( \xi \right) \leq \min_{\xi \in \Xi _{b}}\text{\textsc{var}}%
\left( \xi \right) \leq \text{\textsc{var}}\left( \xi _{0}\right) ,
\label{string}
\end{equation}%
so that we must have equality throughout in (\ref{string}) and $\xi _{0}$ is
\textsc{rbb}($b^{2}$). \hfill $\square $\medskip

\begin{lemma}
\label{lemma: BV}(a) $\xi _{\nu }$ is a solution to (S) for $s^{2}\left( \nu
\right) $. (b) If $s^{2}\geq s^{2}\left( 1\right) =$ \textsc{var}$\left( \xi
_{1}\right) $ then $\xi _{1}$ is a solution to (S) for $s^{2}$.\noindent
\end{lemma}

\noindent The proof of Lemma \ref{lemma: BV} is identical to that of Lemma %
\ref{lemma: BB}, apart from the obvious interchanges \textsc{maxbias }$%
\leftrightarrow $ \textsc{var}, $b\leftrightarrow s$, $\xi
_{0}\leftrightarrow \xi _{1}$, and so is omitted. \medskip

\noindent \textbf{Remark}: The solutions in Theorem \ref{thm: solution}
suggest that the maxima of $b^{2}\left( \cdot \right) $ and $s^{2}\left(
\cdot \right) $ are unique, with%
\begin{equation}
\text{(i) }\{0\}=\arg \max_{\nu \in \left[ 0,1\right] }b^{2}\left( \nu
\right) \text{ and (ii) }\{1\}=\arg \max_{\nu \in \left[ 0,1\right]
}s^{2}\left( \nu \right) .  \label{counter}
\end{equation}%
If not, for instance if $b^{2}\left( \nu \right) $ has multiple maxima or if
`$0$' is not a maximum, then on the set\linebreak $N_{b}=\{ \nu ^{\ast }\in \left[
0,1\right] \left\vert b^{2}\left( \nu ^{\ast }\right) \geq b^{2}\left(
0\right) \right.\} $, both $\xi _{\nu ^{\ast }}$ and $\xi _{0}$ are 
\textsc{rbb}$\left( b^{2}\left( \nu ^{\ast }\right) \right) $, by (a) of
Theorem \ref{thm: solution} and Lemma \ref{lemma: special} respectively,
hence furnish the same minimum variance. Similarly, if (ii) fails then the
bias is constant on the corresponding set $N_{s}$, where $s^{2}\left( \nu
^{\ast }\right) \geq s^{2}\left( 1\right) $. While counterintuitive if these
sets are not singletons -- and if the design weights $\xi _{i}$ vary
continuously with $\nu $ -- these events cannot be ruled out without further
restrictions. This is shown by the example of regression through the origin
-- $p=q=1$, $f\left( x\right) =x$, and $\Xi $ the class of designs placing
mass $\left\{ \alpha ,1-2\alpha ,\alpha \right\} $ on the points of $\chi
=\left\{ -1,0,1\right\} $. For any such design we find that \textsc{var}$%
\left( \xi \right) =1/\alpha $ and \textsc{maxbias}$\left( \xi \right) =1$,
so that $I_{\nu }\left( \xi \right) =\left( 1-\nu \right) /\alpha +\nu $ is
minimized by $\xi _{\nu }=.5\delta _{\pm 1}$ for any $\nu $. Thus $I_{\nu
}\left( \xi _{\nu }\right) =2-\nu $, (\ref{counter}) fails and in fact $%
b^{2}\left( \cdot \right) $ and $s^{2}\left( \cdot \right) $ are constant on 
$\left[ 0,1\right] $: $b^{2}\left( \nu \right) =$ \textsc{maxbias}$\left(
\xi _{\nu }\right) \equiv 1,$ $s^{2}\left( \nu \right) =$ \textsc{var}$%
\left( \xi _{\nu }\right) \equiv 2$.

\section{Examples \label{section: examples}}

We now assume that the design space is finite: $\mathcal{X}=\left\{ 
\boldsymbol{x}_{1},\cdot \cdot \cdot ,\boldsymbol{x}_{N}\right\} $. Let $%
\boldsymbol{Q}$ be an $N\times p$ matrix whose orthonormal columns span the
column space of $\left( \boldsymbol{f}\left( \boldsymbol{x}%
_{1}\right) ,\cdot \cdot \cdot ,\boldsymbol{f}\left( \boldsymbol{x}%
_{N}\right) \right) ^{\prime }$, assumed to have dimension $p$. For a design $\xi 
$ placing mass $\xi _{i}$ on $\boldsymbol{x}_{i}$ define $\boldsymbol{D}%
\left( \xi \right) =diag\left( \xi _{1},\cdot \cdot \cdot ,\xi _{N}\right) $%
. Then in terms of 
\begin{equation*}
\boldsymbol{R}\left( \xi \right) =\boldsymbol{Q}^{\prime }\boldsymbol{D}%
\left( \xi \right) \boldsymbol{Q},\text{ }\boldsymbol{S}\left( \xi \right) =%
\boldsymbol{Q}^{\prime }\boldsymbol{D}^{2}\left( \xi \right) \boldsymbol{Q},%
\text{ }\boldsymbol{U}\left( \xi \right) =\boldsymbol{R}^{-1}\left( \xi
\right) \boldsymbol{S}\left( \xi \right) \boldsymbol{R}^{-1}\left( \xi
\right) \text{,}
\end{equation*}%
it is shown in \cite{w} that 
\begin{equation*}
\text{\textsc{var}}\left( \xi \right) =tr\boldsymbol{R}^{-1}\left( \xi
\right) ,\text{ \textsc{maxbias}}\left( \xi \right) =ch_{\max }\boldsymbol{U}%
\left( \xi \right) .
\end{equation*}%
Here $tr$ and $ch_{\max }$ denote the trace and maximum eigenvalue,
respectively.

The minimization of $I_{\nu }\left( \xi \right) $ is carried out
sequentially, as described in Theorem 5 of  \cite{w}. Briefly, given a current $n$%
-point design $\xi _{n}$, the loss resulting from the addition of a design
point at $\boldsymbol{x}_{i}$ is expanded as 
\begin{equation}
I_{\nu }\left( \xi _{n+1}^{(i)}\right) =I_{\nu }\left( \xi _{n}\right)
-t_{n,i}/n+O\left( n^{-2}\right) ,  \label{expansion}
\end{equation}%
and then $\boldsymbol{x}_{(i)}$, with $\left( i\right) =\arg \max_{i}t_{n,i}$%
, is added to the design. This is carried out to convergence.

A useful measure in choosing a value of the bias/variance parameter $\nu $ is the dimensionless \textit{coefficient of maximum bias}, defined by 
\begin{equation*}
\textsc{cmb}\left( \nu \right) =\sqrt {b^2\left ( \nu \right )/s^2\left( \nu \right)}.
\end{equation*}%
This is akin to the commonly used \textit{coefficient of variation}, but indexes the more relevant -- in our applications -- worst case of the bias.

\begin{figure}
\centering
\includegraphics[scale=1]{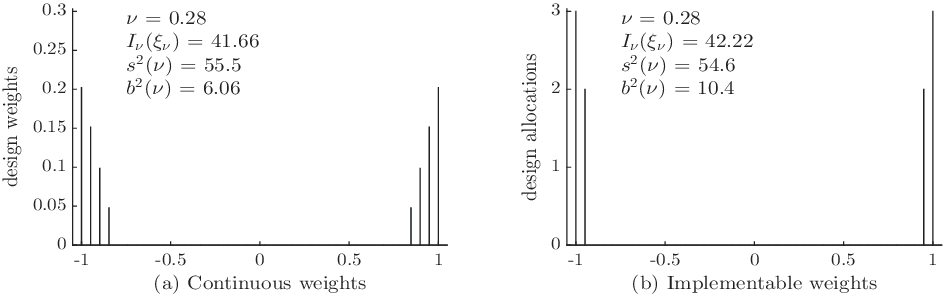}
\caption{(a)\textsc{\ rbb}($b^{2}(.28)$) (also \textsc{rbv}($s^{2}(.28)$)) designs
for the values displayed (\textsc{cmb }= .33). (b) Implementation of the
design in (a); design size $n=10$. }
\label{fig:designs}
\end{figure}

\begin{figure}
\centering
\includegraphics[scale=1]{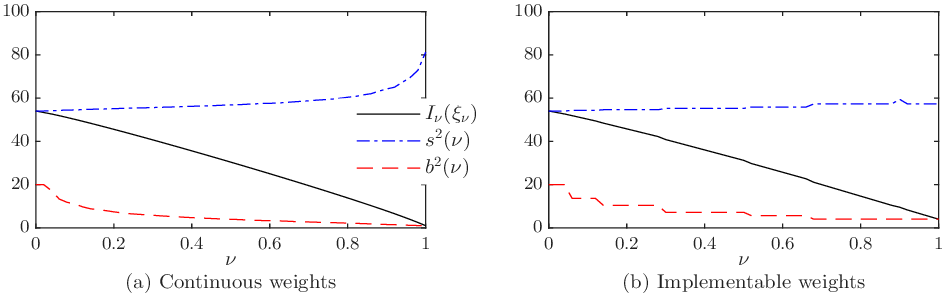}
\caption{\textsc{imse},\textsc{\ var} and \textsc{maxbias} vs.\ $\protect\nu 
$ for the continuous optimal designs and their implementable approximations (%
$n=10$).}
\label{fig:losses}
\end{figure}

\begin{figure}
\centering
\includegraphics[scale=1]{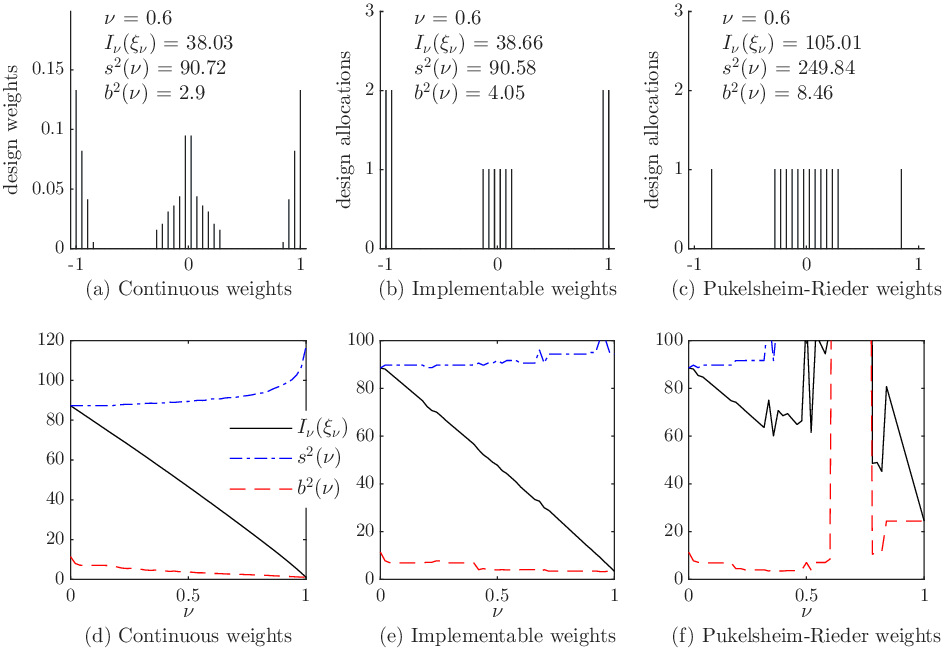}
\caption{(a) -- (c): Continuous and implementable designs for quadratic regression; $\
n=14$. (d) -- (f) \textsc{imse},\textsc{\ var} and \textsc{maxbias} vs.\ $\protect\nu 
$ }
\label{fig:3designs}
\end{figure}

We first present output for the case of straight line regression on a
symmetric design space of size $N=40$ in $\left[ -1,1\right] $: $\chi
=\left\{ x_{i}\left\vert i=1,...,N\right. \right\} $ with $%
x_{i}=-1+2(i-1)/\left( N-1\right) =$ $-x_{N-i+1}$. Put $\boldsymbol{x}%
=\left( x_{1},...,x_{N}\right) ^{\prime }$. Then $\boldsymbol{Q}_{N\times
2}=( \boldsymbol{1}_{N}/\sqrt{N}\text{ }\vdots \text{ }\boldsymbol{x}/%
\sqrt{\left\Vert \boldsymbol{x}\right\Vert }) $. To obtain the \textsc{%
rbb} and \textsc{rbv} designs we minimize $I_{\nu }\left( \xi \right) $ over
unrestricted designs $\xi $, so that the $\xi _{i}$ vary freely over the
simplex $\mathcal{S}=\left\{ \xi _{i}\in \left[ 0,1\right] ,\text{ }\sum \xi
_{i}=1\right\} $. In particular, the minimax design when $\nu =1$ is uniform: $%
\xi _{1,i}\equiv 1/N$, so that $\boldsymbol{U}\left( \xi _{1}\right) =%
\boldsymbol{I}_{N}$ with \textsc{bias}$\left( \xi _{1}\right) =1$ and 
\textsc{var}$\left( \xi _{1}\right) =2N$. When $\nu =0$, the minimax design
is $\xi _{0}=.5\delta _{\pm 1}$, with 
\begin{equation*}
\text{\textsc{maxbias}}\left( \xi _{0}\right) =\max \left( N/2,\left\Vert 
\boldsymbol{x}\right\Vert ^{2}/2\right) =N/2,\quad \text{\textsc{var}}\left(
\xi _{0}\right) =N+\left\Vert \boldsymbol{x}\right\Vert ^{2}.
\end{equation*}%
Thus $\xi _{1}$ is \textsc{rbb}($b^{2}=1$) and $\xi _{0}$ is \textsc{rbv}($%
s^{2}=N+\left\Vert \boldsymbol{x}\right\Vert ^{2}$). Representative results
are presented in parts (a) of Figure \ref{fig:designs} (for which we
specified \textsc{cmb}$\left( \nu \right) \cong 1/3$ and obtained $\nu =.28$%
) and Figure \ref{fig:losses} -- note that (\ref{counter}) holds.

The design in (a) of Figure \ref{fig:designs} is not implementable since the
allocations $n_{i}=$ $n\xi _{i}$ need not be integers. To obtain the implementation in part (b)
of this figure we first rounded up the $n_{i}$ to $\left\lceil n\xi
_{i}\right\rceil $, whose sum then exceeds $n$. The excess is decreased
stepwise, by removing points whose value of $t_{n,i}$ in (\ref{expansion})
is a minimum. This method typically results in only a very small increase in
the minimized value of the \textsc{imse}.

The behaviour shown in (b) of Figure \ref{fig:losses}, in particular of the
bias and as anticipated in the Remark of \S \ref{section: bounded}, reflects
the lack of continuity of the allocations as functions of $\nu $. Together
with \textsc{cmb}$\left( \cdot \right) $ this plot can also serve as a guide
to the designer in choosing a value of $\nu $ for an implementable design.
\medskip

\noindent \textbf{Remark}: Our method of rounding the design weights so as
to obtain implementable designs is somewhat non-standard, and is intended to
preserve, as much as possible, the minimized \textsc{imse}. A more common
method is the `efficient design apportionment' method of \cite{pr92}. This
is a rounding procedure that has, amongst others, the property of `sample
size monotonicity' -- if a new point is to be added to an existing design,
then none of the current allocations will be reduced. Unless this property
is required we cannot recommend this method in the current application, as
it is too often very unstable, resulting in large increases in the loss. An
example is quadratic regression on the same design space as above, with a
design size $n=14$, as illustrated in Figure \ref{fig:3designs}.

\section*{Acknowledgements}

This work was carried out with the support of the Natural Sciences and Engineering Council of Canada. It has benefited from the helpful comments of two anonymous referees.

\bibliographystyle{natbib}
\bibliography{references}
\end{document}